\documentclass[a4paper,10pt]{amsart}

\usepackage{amssymb}
\usepackage{fullpage}
\usepackage[OT4]{fontenc}
\usepackage{texdraw}

\theoremstyle{definition}
\newtheorem{definition}{Definition}
\theoremstyle{plain}
\newtheorem{theorem}[definition]{Theorem}
\newtheorem{proposition}[definition]{Proposition}

\newtheorem{conjecture}[definition]{Conjecture}

\theoremstyle{remark}

\newtheorem{example}[definition]{Example}

\DeclareMathOperator{\edim}{edim}
\DeclareMathOperator{\vdim}{vdim}

\DeclareMathOperator{\mult}{mult}

\DeclareMathOperator{\red}{red}
\DeclareMathOperator{\Pic}{Pic}
\DeclareMathOperator{\Cr}{Cr}
\DeclareMathOperator{\CrSt}{Cr^{\circ}}

\def\field{\mathbb{K}}
\def\N{\mathbb{N}}
\def\TT{\mathbb{T}}
\def\Z{\mathbb{Z}}

\def\PP{\mathbb{P}}

\def\sys{\mathcal{L}}

\def\rdf{:=}
\def\dff{\it}

\let\to\longrightarrow
\let\mapsto\longmapsto

\def\epar{\smallskip}

\newenvironment{algorithm}[1]{\medskip\noindent{\bf Algorithm} #1\\}{\medskip}

\begin{document}

\title{Quasi-homogeneous linear systems on $\mathbb P^2$ with base
points of multiplicity 7, 8, 9, 10}

\author{Marcin Dumnicki}

\dedicatory{
Institute of Mathematics, Jagiellonian University, \\
Reymonta 4, 30-059 Krak\'ow, Poland \\
}

\thanks{Email address: Marcin.Dumnicki@im.uj.edu.pl}





\begin{abstract}
In the paper we prove Harbourne-Hirschowitz conjecture
for quasi-homogeneous linear systems on $\mathbb P^2$ for
$m=7$, 8, 9, 10, i.e.
systems of curves of given degree passing through points in general
position with multiplicities at least $m,\dots,m,m_0$, where $m=7$, 8, 9, 10, $m_0$
is arbitrary.
\end{abstract}

\maketitle

\section{Introduction}

In what follows we assume that the ground field $\field$ is of characteristic zero.
Let $d \in \Z$, let $m_1,\dots,m_r \in \N$.
By $\sys(d;m_1,\dots,m_r)$ we denote the linear system of curves
(in $\PP^2 \rdf \PP^2(\field)$)
of degree $d$ passing through $r$ points $p_1,\dots,p_r$ in general position with
multiplicities at least $m_1,\dots,m_r$. The dimension of such system is denoted by
$\dim \sys(d;m_1,\dots,m_r)$.
Define the {\dff virtual dimension of $L = \sys(d;m_1,\dots,m_r)$}
\begin{align*}
\vdim L & \rdf \binom{d+2}{2} - \sum_{j=1}^{r} \binom{m_j+1}{2} - 1 \\
\intertext{and the {\dff expected dimension of $L$}}
\edim L & \rdf \max \{ \vdim L, -1 \}.
\end{align*}
Observe that $\dim L \geq \edim L$. If this inequality is strict then
$L$ is called {\dff special}, {\dff non-special} otherwise.
The system $L$ is called {\dff non-empty} if $\dim L \geq 0$, {\dff empty} otherwise.

Let $\pi:X \to \PP^2$ be the blow-up of $\PP^2$ at $r$ points in general position.
The Picard group $\Pic(X)$ of $X$ is generated by $H,E_1,\dots,E_r$, where
$H$ is the pullback of the class of a line in $\PP^2$, $E_1,\dots,E_r$
are exceptional divisors. The system $\sys(d;m_1,\dots,m_r)$
is isomorphic to the complete linear system (on $X$) associated to the divisor
$dH-m_1E_1-\dots-m_rE_r$. Observe that this way we can define
$\sys(d;m_1,\dots,m_r)$ for $m_1,\dots,m_r \in \Z$.
Consider the standard intersection form on $X$ given by
$H^2=1$, $E_j^2=-1$, $H.E_j=0$, $E_j.E_\ell=0$ for $j \neq \ell$.
Now (by Riemann-Roch)
$$\vdim L = \frac{L^2-L.K_{X}}{2},$$
where $K_{X}$ is the canonical divisor on $X$.
In what follows we always allow (unless stated otherwise) 
negative multiplicities.

We say that a curve $E \subset X$ is
a {\dff $-1$-curve} on $X$ if $E$ is irreducible and $E^2=E.K_{X}=-1$.
We recall the following definition of $-1$-special system (see e.g. \cite{CMirdeg}):

\begin{definition}
A linear system $L$ is {\dff $-1$-special} if there exists
$-1$-curves $C_1,\dots,C_s$ such that $L.C_j = -k_j$, $k_j \geq 1$ for $j=1,\dots,s$,
$k_j \geq 2$ for some $j$ and the system
$M=L-(k_1C_1+\dots+k_sC_s)$ has non-negative virtual dimension and
non-negative intersection with every $-1$-curve.
\end{definition}

From the above definition it is clear that if $L$ is $-1$-special then it is
non-empty and its dimension is at least $\dim M$. Since computation of
$\vdim M$ leads to the inequality $\vdim M > \vdim L$, every $-1$-special
system is special. Being more precise, if $L.C_j=-k_j$, for $j=1,\dots,s$, then
$$\vdim L = \vdim M + \sum_{j=1}^{s} \frac{k_j-k_j^2}{2}.$$

The converse is only conjectured to hold:

\begin{conjecture}[Harbourne-Hirschowitz]
A linear system is special if and only if it is $-1$-special.
\end{conjecture}

The above conjecture is known to hold in some cases.
The case $r \leq 9$ has been solved by Nagata (\cite{Nag2}).
For low multiplicities (i.e. bounded by some constant)
it begun with \cite{HirF}, where the case
$m_1=\dots=m_r \leq 3$ was solved. 
The case when all multiplicities are bounded by $4$ has been done by
\cite{Mignon}, it has been extended to $7$ in \cite{Yang}
and $11$ in \cite{mdwj}.

The homogeneous case ($m_1=\dots=m_r$)
with multiplicities up to $42$ has been succesfully solved in \cite{md}.
The quasi-homogeneous case ($m_1=\dots=m_{r-1}$, $m_r$ arbitrary) has been
done 
for $m_1=3$ in \cite{CMirdeg}, for $m_1=4$ in \cite{sei4,laf4},
for $m_1=5$ in \cite{base5}, for $m_1=6$ in \cite{new6}. Our result is the following:

\begin{theorem}
\label{mainth}
The Harbourne-Hirschowitz conjecture holds for quasi-homogeneous systems
with almost all multiplicities equal to $7$, $8$, $9$ or $10$.
\end{theorem}

The methods showing the conjecture for $m_1=3,4,5,6$ used by
authors mentioned above are of the same type.
Namely, using degeneration method introduced by Ciliberto and Miranda
we can show non-specialty of a large family of systems
with many base points. In fact, with the help of this method we can show that if 
the family of systems 
$$\{\sys(d;m^{k},m_0) : d,m_0 \in \N, k_1\leq k \leq k_2\},$$
for a carefully chosen $k_1$ and $k_2$, contains only non-special ones
then all systems of the form
$\sys(d;m^{\times \ell},m_0)$ for $\ell \geq k_2$ are non-special.
Another task is to show that if the difference between $d$ and $m_0$ in
a system $\sys(d;m^{\times k},m_0)$ is big enough then the system is
non-special.

Having shown the above, we are left with a family of cases that
can be solved using degeneration method, Cremona transformation,
``ad hoc'' arguments and computations of the rank of an interpolation matrices.

Authors of \cite{CMirdeg,sei4,laf4,base5,new6} used computer programs to deal with
large number of cases. The programs are of two kinds.
The first one, it is an implementation of degeneration technique --- for large
number of cases we must check whether degeneration exists or not.
The result (for a single case) can be easily checked by hand, the reason
for using software is the number of cases.
The second kind uses computer programs to evaluate the dimension of
a given system by a direct computation of the rank of an appropriate matrix,
which, in interesting cases, has large size (e.g. $105 \times 105$).
This cannot be done by hand, from obvious reasons.

In this paper we use the same approach, but another methods.
Instead of degeneration method we will use reduction algorithm
introduced in \cite{redmd} and \cite{mdwj} together with
direct computations of dimension of systems.
To deal with the remaining cases we will use Cremona transformation
(see Def. \ref{defcrem}), ``glueing'' of points (see Thm. \ref{glue}) and known results.

We note here that both approaches, by degeneration and reduction algorithm,
promise to be usable for larger values of $m$ (quasi-homogeneous multiplicity).
We prefer the second one --- observe that this paper is not much longer than
\cite{sei4,base5,new6}, however, we deal with four bigger multiplicities at once.

The paper is organized as follows: The next section is devoted
to present some methods of showing $-1$-specialty or non-specialty of systems.
Section 3 contains a brief introduction to the reduction method together
with the results obtained with the help of this method and computer programs.
In section 4 we deal with the remaining cases, i.e. systems with
few base points and low difference between the degree and quasi-homogeneous
multiplicity.
The last section contains a note on Seibert's work.

\section{Tools}

\begin{theorem}[splitting]
\label{maintool}
Let $d,k,m_1,\dots,m_r,m_{r+1},\dots,m_s \in \N$.
If
\begin{itemize}
\item
$L_1 = \sys(k;m_{1},\dots,m_s)$ is non-special,
\item
$L_2 = \sys(d;m_{s+1},\dots,m_r,k+1)$ is non-special,
\item
$(\vdim L_1+1)(\vdim L_2+1) \geq 0$,
\end{itemize}
then the system $L=\sys(d;m_1,\dots,m_r)$ is non-special.
\end{theorem}

\begin{proof}
See \cite{mdprep}, Thm. 1.
\end{proof}

\begin{theorem}[glueing]
\label{glue}
Let $\sys(k;m^s)$ be non-special,
let 
\begin{align*}
L_1 & = \sys(d;m_1,\dots,m_r,m^s),\\
L_2 & = \sys(d;m_1,\dots,m_r,k+1).
\end{align*}
If either $-1 \leq \vdim L_2 \leq \vdim L_1$ or $\vdim L_1 \leq \vdim L_2 \leq -1$
then in order to show non-specialty of $L_1$ it is enough to show non-specialty
of $L_2$.
\end{theorem}

\begin{proof}
Follows from Theorem \ref{maintool}.
\end{proof}

\begin{definition}
\label{defcrem}
Let $d,m_1,\dots,m_r \in \Z$, let
$k = d-(m_1+m_2+m_3)$.
Define the {\dff Cremona transformation of a system $\sys(d;m_1,\dots,m_r)$}
$$\Cr(\sys(d;m_1,\dots,m_r)) \rdf \sys(d+k;m_1+k,m_2+k,m_3+k,m_4,\dots,m_r).$$
\end{definition}

\begin{theorem}
\label{cremona}
Let $L$ be a linear system. The following holds:
\begin{enumerate}
\item
$\dim \Cr(L) = \dim L$,
\item
$L$ is special if and only if $\Cr(L)$ is special,
\item
$L$ is $-1$-special if and only if $\Cr(L)$ is $-1$-special,
\end{enumerate}
\end{theorem}

\begin{proof}
The proof can be found, for example, in \cite{thin}.
The idea is to show that standard birational transformation induces
an action on $\Pic(X)$ such that
$H \mapsto 2H-E_1-E_2-E_3$, $E_j \mapsto H-(E_1+E_2+E_3)+E_j$ for $j=1,\dots,3$
and $E_j \mapsto E_j$ for $j \geq 4$. Observe that $-1$-curves are
transformed into $-1$-ones.
\end{proof}

\begin{definition}
We say that {\dff $\sys_n(d;m_1,\dots,m_r)$ is in standard form}
if $d < 0$ or the following holds:
\begin{itemize}
\item
$m_1,\dots,m_r$ are non-increasing,
\item
$d - (m_1+m_2+m_3) \geq 0$.
\end{itemize}
\end{definition}

Every system can be transformed (by a finite number of
Cremona transformations and sorting of multiplicities)
into a standard form. For a system $L$ we denote its standard form
by $\CrSt(L)$.

Let $L=\sys(d;m_1,\dots,m_r)$ be a linear system in standard form. 
From \cite{thin} we may understand
what happens if some of $d,m_1,\dots,m_r$ are negative.
\begin{enumerate}
\item
if $d<0$ then $L$ is empty;
\item
if $m_j = -1$ then $E_j$ is a fixed component for $L$; 
let $L'=\sys(d;m_1,\dots,m_{j-1},0,m_{j+1},\dots,m_r)$.
Since $\vdim L=\vdim L'$, $\dim L=\dim L'$, it is enough to
study $L'$;
\item
if $m_j \leq -2$ then $E_j$ is a multiple fixed component; since $E_j$ is
an $-1$-curve, the system $L$ is special  
if and only if $L'=\sys(d;m_1,\dots,m_{j-1},0,m_{j+1},\dots,m_r)$ is non-empty.
Moreover, if $L'$ is non-empty and non-special, or it is $-1$-special then $L$
is $-1$ special.
\end{enumerate}

Moreover, (see \cite{thin})
the intersection number $L.C$, where $C$ is an $-1$-curve, is non-negative for
any system in standard form with non-negative multiplicities, hence
such system cannot be $-1$-special. If, additionally, it is a system
for which the Harbourne-Hirschowitz conjecture has been proved
(e.g. multiplicities bounded by $11$ or based on at most $9$ points) then
it is non-special.

We recall the following result (which has been mentioned in the introduction).

\begin{theorem}[\cite{mdwj}]
\label{up11}
The Harbourne-Hirschowitz conjecture is true for
systems with multiplicities bounded by $11$.
\end{theorem}

\section{Results using reduction method}

The first step is to show that systems with large number of points are
non-special. To do this we will use reduction method introduced
in \cite{redmd} and then exploited in \cite{mdwj}.

For a finite $D \subset \TT^2 \rdf \{ x^{\alpha_1}y^{\alpha_2} \subset \field[x,y] : \alpha_1,\alpha_2 \in \N \}$
and multiplicities $m_1,\dots,m_r$ define the space
$$V(D;m_1,\dots,m_r) \rdf \bigg\{ f=\sum_{t \in D} c_tt \in \field[x,y] : \mult_{p_j}f \geq m_j\bigg\}$$
for $p_1,\dots,p_r \in \field^2$ in general position.
We say that $V(D;m_1,\dots,m_r)$ is {\dff non-special} if its dimension
(as a vector space over $\field$) is equal to its expected dimension
\begin{align*}
\edim V(D;m_1,\dots,m_r) & \rdf \max \{\vdim V(D;m_1,\dots,m_r), 0 \}, \\
\vdim V(D;m_1,\dots,m_r) & \rdf \#D-\sum_{j=1}^r \binom{m_j+1}{2}.
\end{align*}
Observe that $\sys(d;m_1,\dots,m_r)$ is non-special if and only if
$V(D;m_1,\dots,m_r)$ is non-special for $D=\{t \in \TT^2 : \deg t \leq d\}$.

\begin{definition}
Let $a_{1}, \dots, a_{k} \in \N$, $a_{j} \leq j$, $j=1,\dots,k$.
Define the {\dff diagram $(a_{1},\dots,a_{k})$}
$$(a_{1},\dots,a_{k}) = \bigcup_{j=1}^{k} \{ x^{\alpha_{1}}y^{\alpha_{2}} 
\in \TT^{2} : \alpha_{1} + \alpha_{2} = j - 1, \
\alpha_{2} < a_{j} \}.$$
A single set
$\{ x^{\alpha_{1}}y^{\alpha_{2}} 
\in \TT^{2} : \alpha_{1} + \alpha_{2} = j - 1, \
\alpha_{2} < a_{j} \}$
will be called a {\dff $j$-th layer}, or simply a {\dff layer}.
For $a, a_{1}, \dots, a_{k} \in \N$, $a_{j} \leq a+j$ define
$$(\overline{a},a_{1},\dots,a_{k}) := (1,2,\dots,a-1,a,a_{1},\dots,a_{k}).$$
\end{definition}

We will also use notation
$$(\overline{a},\{b\}^{\times p},a_1,a_2,\dots) \rdf
(\overline{a},\underbrace{b,\dots,b}_{p \text{ times}},a_1,a_2,\dots).$$
Observe that for $d \geq 1$ we have
$\{t \in \TT^2 : \deg t \leq d\} = (\overline{d+1})$.

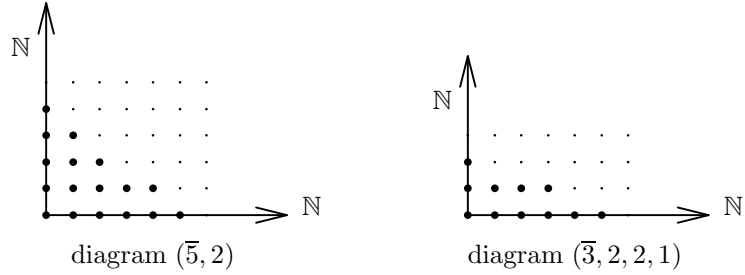
\begin{example}
$$
\begin{array}{ccc}
\begin{texdraw}
\drawdim pt
\move(0 0)
\fcir f:0 r:0.5
\move(0 10)
\fcir f:0 r:0.5
\move(0 20)
\fcir f:0 r:0.5
\move(0 30)
\fcir f:0 r:0.5
\move(0 40)
\fcir f:0 r:0.5
\move(0 50)
\fcir f:0 r:0.5
\move(10 0)
\fcir f:0 r:0.5
\move(10 10)
\fcir f:0 r:0.5
\move(10 20)
\fcir f:0 r:0.5
\move(10 30)
\fcir f:0 r:0.5
\move(10 40)
\fcir f:0 r:0.5
\move(10 50)
\fcir f:0 r:0.5
\move(20 0)
\fcir f:0 r:0.5
\move(20 10)
\fcir f:0 r:0.5
\move(20 20)
\fcir f:0 r:0.5
\move(20 30)
\fcir f:0 r:0.5
\move(20 40)
\fcir f:0 r:0.5
\move(20 50)
\fcir f:0 r:0.5
\move(30 0)
\fcir f:0 r:0.5
\move(30 10)
\fcir f:0 r:0.5
\move(30 20)
\fcir f:0 r:0.5
\move(30 30)
\fcir f:0 r:0.5
\move(30 40)
\fcir f:0 r:0.5
\move(30 50)
\fcir f:0 r:0.5
\move(40 0)
\fcir f:0 r:0.5
\move(40 10)
\fcir f:0 r:0.5
\move(40 20)
\fcir f:0 r:0.5
\move(40 30)
\fcir f:0 r:0.5
\move(40 40)
\fcir f:0 r:0.5
\move(40 50)
\fcir f:0 r:0.5
\move(50 0)
\fcir f:0 r:0.5
\move(50 10)
\fcir f:0 r:0.5
\move(50 20)
\fcir f:0 r:0.5
\move(50 30)
\fcir f:0 r:0.5
\move(50 40)
\fcir f:0 r:0.5
\move(50 50)
\fcir f:0 r:0.5
\move(60 0)
\fcir f:0 r:0.5
\move(60 10)
\fcir f:0 r:0.5
\move(60 20)
\fcir f:0 r:0.5
\move(60 30)
\fcir f:0 r:0.5
\move(60 40)
\fcir f:0 r:0.5
\move(60 50)
\fcir f:0 r:0.5
\arrowheadtype t:V
\move(0 0)
\avec(90 0)
\move(0 0)
\avec(0 80)
\htext(96 0){$\mathbb{N}$}
\htext(-13 60){$\mathbb{N}$}
\move(0 0)
\fcir f:0 r:1.5
\move(10 0)
\fcir f:0 r:1.5
\move(0 10)
\fcir f:0 r:1.5
\move(20 0)
\fcir f:0 r:1.5
\move(10 10)
\fcir f:0 r:1.5
\move(0 20)
\fcir f:0 r:1.5
\move(30 0)
\fcir f:0 r:1.5
\move(20 10)
\fcir f:0 r:1.5
\move(10 20)
\fcir f:0 r:1.5
\move(0 30)
\fcir f:0 r:1.5
\move(40 0)
\fcir f:0 r:1.5
\move(30 10)
\fcir f:0 r:1.5
\move(20 20)
\fcir f:0 r:1.5
\move(10 30)
\fcir f:0 r:1.5
\move(0 40)
\fcir f:0 r:1.5
\move(50 0)
\fcir f:0 r:1.5
\move(40 10)
\fcir f:0 r:1.5
\textref h:C v:C
\htext(40 -15){diagram $(\overline{5},2)$}
\end{texdraw}
&
\hspace{1cm}
&
\begin{texdraw}
\drawdim pt
\move(0 0)
\fcir f:0 r:0.5
\move(0 10)
\fcir f:0 r:0.5
\move(0 20)
\fcir f:0 r:0.5
\move(0 30)
\fcir f:0 r:0.5
\move(10 0)
\fcir f:0 r:0.5
\move(10 10)
\fcir f:0 r:0.5
\move(10 20)
\fcir f:0 r:0.5
\move(10 30)
\fcir f:0 r:0.5
\move(20 0)
\fcir f:0 r:0.5
\move(20 10)
\fcir f:0 r:0.5
\move(20 20)
\fcir f:0 r:0.5
\move(20 30)
\fcir f:0 r:0.5
\move(30 0)
\fcir f:0 r:0.5
\move(30 10)
\fcir f:0 r:0.5
\move(30 20)
\fcir f:0 r:0.5
\move(30 30)
\fcir f:0 r:0.5
\move(40 0)
\fcir f:0 r:0.5
\move(40 10)
\fcir f:0 r:0.5
\move(40 20)
\fcir f:0 r:0.5
\move(40 30)
\fcir f:0 r:0.5
\move(50 0)
\fcir f:0 r:0.5
\move(50 10)
\fcir f:0 r:0.5
\move(50 20)
\fcir f:0 r:0.5
\move(50 30)
\fcir f:0 r:0.5
\move(60 0)
\fcir f:0 r:0.5
\move(60 10)
\fcir f:0 r:0.5
\move(60 20)
\fcir f:0 r:0.5
\move(60 30)
\fcir f:0 r:0.5
\arrowheadtype t:V
\move(0 0)
\avec(90 0)
\move(0 0)
\avec(0 60)
\htext(96 0){$\mathbb{N}$}
\htext(-13 40){$\mathbb{N}$}
\move(0 0)
\fcir f:0 r:1.5
\move(10 0)
\fcir f:0 r:1.5
\move(0 10)
\fcir f:0 r:1.5
\move(20 0)
\fcir f:0 r:1.5
\move(10 10)
\fcir f:0 r:1.5
\move(0 20)
\fcir f:0 r:1.5
\move(30 0)
\fcir f:0 r:1.5
\move(20 10)
\fcir f:0 r:1.5
\move(40 0)
\fcir f:0 r:1.5
\move(30 10)
\fcir f:0 r:1.5
\move(50 0)
\fcir f:0 r:1.5
\textref h:C v:C
\htext(40 -15){diagram $(\overline{3},2,2,1)$}
\end{texdraw}
\end{array}
$$
\end{example}

\begin{definition}
\label{defred}
Let $m \in \N^{*}$, let 
$D = (b_{1},\dots,b_{\ell},a_{1},\dots,a_{m})$ be a diagram, $a_{m} > 0$.
Define the numbers $v_{j} \in \N$, $j=1,\dots,m$ together with sets
$V_{j}$, $j=0,\dots,m$ inductively (beginning with $m$, going down to $0$) to be:
\begin{align*}
V_{m} & := \{1,\dots,m\}, \\
v_{j} & := \left\{
\begin{array}{ll}
a_{j}, & a_{j} < m \textrm{ and } \max V_{j} \geq a_{j} \\
\max V_{j},& \textrm{otherwise.}
\end{array}
\right. \\
V_{j-1} & := V_{j} \setminus \{v_{j}\}.\\
\end{align*}
If we have $V_{0} = \varnothing$
then we say that
\textit{$D$ is $m$-reducible}. The diagram
$$\red_{m}(D) := (b_{1},\dots,b_{\ell},
a_{1}-v_{1}, \dots, a_{m}-v_{m})$$
will be called the \textit{$m$-reduction of $D$}. 
\end{definition}

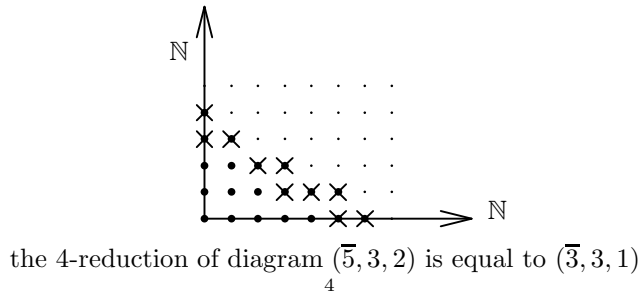
\begin{example}
$$
\begin{texdraw}
\drawdim pt
\move(0 0)
\fcir f:0 r:0.5
\move(0 10)
\fcir f:0 r:0.5
\move(0 20)
\fcir f:0 r:0.5
\move(0 30)
\fcir f:0 r:0.5
\move(0 40)
\fcir f:0 r:0.5
\move(0 50)
\fcir f:0 r:0.5
\move(10 0)
\fcir f:0 r:0.5
\move(10 10)
\fcir f:0 r:0.5
\move(10 20)
\fcir f:0 r:0.5
\move(10 30)
\fcir f:0 r:0.5
\move(10 40)
\fcir f:0 r:0.5
\move(10 50)
\fcir f:0 r:0.5
\move(20 0)
\fcir f:0 r:0.5
\move(20 10)
\fcir f:0 r:0.5
\move(20 20)
\fcir f:0 r:0.5
\move(20 30)
\fcir f:0 r:0.5
\move(20 40)
\fcir f:0 r:0.5
\move(20 50)
\fcir f:0 r:0.5
\move(30 0)
\fcir f:0 r:0.5
\move(30 10)
\fcir f:0 r:0.5
\move(30 20)
\fcir f:0 r:0.5
\move(30 30)
\fcir f:0 r:0.5
\move(30 40)
\fcir f:0 r:0.5
\move(30 50)
\fcir f:0 r:0.5
\move(40 0)
\fcir f:0 r:0.5
\move(40 10)
\fcir f:0 r:0.5
\move(40 20)
\fcir f:0 r:0.5
\move(40 30)
\fcir f:0 r:0.5
\move(40 40)
\fcir f:0 r:0.5
\move(40 50)
\fcir f:0 r:0.5
\move(50 0)
\fcir f:0 r:0.5
\move(50 10)
\fcir f:0 r:0.5
\move(50 20)
\fcir f:0 r:0.5
\move(50 30)
\fcir f:0 r:0.5
\move(50 40)
\fcir f:0 r:0.5
\move(50 50)
\fcir f:0 r:0.5
\move(60 0)
\fcir f:0 r:0.5
\move(60 10)
\fcir f:0 r:0.5
\move(60 20)
\fcir f:0 r:0.5
\move(60 30)
\fcir f:0 r:0.5
\move(60 40)
\fcir f:0 r:0.5
\move(60 50)
\fcir f:0 r:0.5
\move(70 0)
\fcir f:0 r:0.5
\move(70 10)
\fcir f:0 r:0.5
\move(70 20)
\fcir f:0 r:0.5
\move(70 30)
\fcir f:0 r:0.5
\move(70 40)
\fcir f:0 r:0.5
\move(70 50)
\fcir f:0 r:0.5
\arrowheadtype t:V
\move(0 0)
\avec(100 0)
\move(0 0)
\avec(0 80)
\htext(106 0){$\mathbb{N}$}
\htext(-13 60){$\mathbb{N}$}
\move(0 0)
\fcir f:0 r:1.5
\move(10 0)
\fcir f:0 r:1.5
\move(0 10)
\fcir f:0 r:1.5
\move(20 0)
\fcir f:0 r:1.5
\move(10 10)
\fcir f:0 r:1.5
\move(0 20)
\fcir f:0 r:1.5
\move(30 0)
\fcir f:0 r:1.5
\move(20 10)
\fcir f:0 r:1.5
\move(10 20)
\fcir f:0 r:1.5
\move(0 30)
\fcir f:0 r:1.5
\move(40 0)
\fcir f:0 r:1.5
\move(30 10)
\fcir f:0 r:1.5
\move(20 20)
\fcir f:0 r:1.5
\move(10 30)
\fcir f:0 r:1.5
\move(0 40)
\fcir f:0 r:1.5
\move(50 0)
\fcir f:0 r:1.5
\move(40 10)
\fcir f:0 r:1.5
\move(30 20)
\fcir f:0 r:1.5
\move(60 0)
\fcir f:0 r:1.5
\move(50 10)
\fcir f:0 r:1.5
\move(0 30)
\move(-3 27)
\lvec(3 33)
\move(3 27)
\lvec(-3 33)
\move(0 40)
\move(-3 37)
\lvec(3 43)
\move(3 37)
\lvec(-3 43)
\move(10 30)
\move(7 27)
\lvec(13 33)
\move(13 27)
\lvec(7 33)
\move(20 20)
\move(17 17)
\lvec(23 23)
\move(23 17)
\lvec(17 23)
\move(30 10)
\move(27 7)
\lvec(33 13)
\move(33 7)
\lvec(27 13)
\move(30 20)
\move(27 17)
\lvec(33 23)
\move(33 17)
\lvec(27 23)
\move(40 10)
\move(37 7)
\lvec(43 13)
\move(43 7)
\lvec(37 13)
\move(50 0)
\move(47 -3)
\lvec(53 3)
\move(53 -3)
\lvec(47 3)
\move(50 10)
\move(47 7)
\lvec(53 13)
\move(53 7)
\lvec(47 13)
\move(60 0)
\move(57 -3)
\lvec(63 3)
\move(63 -3)
\lvec(57 3)
\textref h:C v:C
\htext(45 -15){the $4$-reduction of diagram $(\overline{5},3,2)$ is equal to $(\overline{3},3,1)$}
\end{texdraw}
$$
As an another example consider the diagram $(\overline{32})$. We will perform
one $12$-reduction and four $9$-reductions. The resulting diagram
is equal to $(\overline{19},18,17,16,14,10,5)$:
$$
\begin{array}{rrrrrrrrrrrrrr}
(\overline{19}, & 20, & 21, & 22, & 23, & 24, & 25, & 26, & 27, & 28, & 29, & 30, & 31, & 32) \\
              &    & -1 & -2 & -3 & -4 & -5 & -6 & -7 & -8 & -9 &-10 &-11 &-12 \\
\hline
(\overline{19}, & 20, & 20, & 20, & 20, & 20, & 20, & 20, & 20, & 20, & 20, & 20, & 20, & 20) \\
              &    &    &    &    & -1 & -2 & -3 & -4 & -5 & -6 & -7 & -8 & -9 \\
\hline
(\overline{19}, & 20, & 20, & 20, & 20, & 19, & 18, & 17, & 16, & 15, & 14, & 13, & 12, & 11) \\
              &    &    &    &    & -1 & -2 & -3 & -4 & -5 & -6 & -7 & -8 & -9 \\
\hline
(\overline{19}, & 20, & 20, & 20, & 20, & 18, & 16, & 14, & 12, & 10, &  8, &  6, &  4, &  2) \\
              &    &    &    &    & -1 & -3 & -5 & -7 & -9 & -8 & -6 & -4 & -2 \\
\hline
(\overline{19}, & 20, & 20, & 20, & 20, & 17, & 13, &  9, &  5, &  1) \\
              & -2 & -3 & -4 & -6 & -7 & -8 & -9 & -5 & -1 \\
\hline
(\overline{19}, & 18, & 17, & 16, & 14, & 10, & 5) \\
\end{array}
$$
We can perform additional five $9$-reductions to obtain $(\overline{6},6,6,5,5,2)$.
The sequence of reductions presented above will be used later to show
non-specialty of $\sys(31;12,9^{\times 9})$.
\end{example}

\begin{definition}
For a $m_r$-reducible diagram $D$ we will say that 
space $V(\red_{m_r}(D);m_1,\dots,m_{r-1})$ is an {\dff $m$-reduction}
of $V(D;m_1,\dots,m_r)$.
\end{definition}

The reduction method is based on the following fact (see \cite{redmd}
for detailed proof; also the sketch of proof can be found in \cite{mdwj}):

\begin{theorem}
\label{usingred}
Let $m_{1},\dots,m_{r} \in \N$. Let $V=V(D;m_{1},\dots,m_{r})$.
If $D$ is $m_{r}$-reducible and the $m_r$-reduction 
of $V$ is non-special then $V$ is non-special.
\end{theorem}

Let $V=V(D;m_1,\dots,m_r)$. We can reduce $V$ until all multiplicities disappear
or the resulting diagram in no longer $m_j$-reducible for all remaining $m_j$'s.
Observe that the $m$-reduction is performed on the last $m$ layers
of a diagram. Let $D=(\dots,b_1,b_2,\dots,b_m)$. We will try to reduce
$D$ to $(\dots,b_1',b_2',\dots,b'_m)$. We have the following:
\begin{itemize}
\item
if $b_j +1 \geq b_{j+1} \geq 2m$ then $m$-reduction on layers
``$b_j$'' and ``$b_{j+1}$'' is possible; after reducing we will have
$b'_j \geq b'_{j+1} \geq m$;
\item
if $b_j \geq b_{j+1} \geq m$ then $m$-reduction on layers
``$b_j$'' and ``$b_{j+1}$'' is possible; after reducing we will have
$b'_j > b'_{j+1}$;
\item
finally, if $b_j > b_{j+1}$ then $m$-reduction on layers
``$b_j$'' and ``$b_{j+1}$'' is possible.
\end{itemize}
In \cite{redmd} one can find additional information on
how long reducing is possible.
We deduce that a diagram $D=(\overline{a},\{a\}^{\times k})$
for $a \geq m$ and given $\ell \leq k$ 
can be reduced (using $m$-reductions; if $k-\ell \leq m-1$ then
we can use no reductions) to a diagram
$(\overline{a},\{a\}^{\times \ell},a_1,\dots,a_{m-1})$.
If, moreover, $a\geq 2m$ and $2m \leq b \leq a$ then $D$ can be reduced to
a diagram $(\overline{b},a_1,\dots,a_{m-1})$.

We will use reductions to show non-specialty of large families of
systems. Before that define (for a diagram $D$ and $m > 0$) the following number:
$$p(D) \rdf \left\lfloor \frac{\# D}{\binom{m+1}{2}} \right\rfloor.$$

\begin{proposition}
\label{downszer}
Let $m > 0$, let $a \geq m$, let $k \geq 0$.
Let 
$$\mathcal D = \{ (\overline{a},\{a\}^{\times k},a_1,\dots,a_{m-1}) : 
a \geq a_1 \geq a_2 \geq \dots \geq a_{m-1} \}.$$
If for all $D \in \mathcal D$ the spaces
$V(D;m^{\times p(D)})$ and $V(D;m^{\times (p(D)+1)})$ are non-special
then for every $r > \max \{p(D) : D \in \mathcal D\}+1$ and
$m_0 \in \N$ the system
$\sys(m_0+a-1;m_0,m^{\times r})$ is non-special.
\end{proposition}

\begin{proof}
We will show that the space $V(D;m_0,m^{\times a-1})$ is non-special
for $D=(\overline{m_0+a})$. We can $m_0$-reduce our space to
$V'=V((\overline{a},\{a\}^{\times \ell});m^{\times r})$ for
some $\ell \geq 0$.
If $\ell < k$ then $\vdim V' < 0$ and, since $r$ is big enough,
the same holds for $V'=V((\overline{a},\{a\}^{\times k});m^{\times r})$.
So, without loss of generality, we may assume $\ell \geq k$.
Performing $m$-reductions on a diagram $(\overline{a},\{a\}^{\times \ell})$
leads to some diagram $D \in \mathcal D$, or we obtain a system
without conditions. In any case, using Thm. \ref{usingred},
we complete the proof.
\end{proof}

\begin{proposition}
\label{downup}
Let $m \in \N$, let $a \geq 2m$.
Let 
$$\mathcal D = \{ (\overline{a},a_1,\dots,a_{m-1}) : 
a \geq a_1 -1 \geq a_2 -2 \geq \dots \geq a_{m-1} -(m-1) \}.$$
If for all $D \in \mathcal D$ the spaces
$V(D;m^{\times p(D)})$ and $V(D;m^{\times (p(D)+1)})$ are non-special
then for every $m_0 \in \N$, $r > \max \{p(D) : D \in \mathcal D\}+1$
and $d\geq m_0+a-1$ the system
$\sys(d;m_0,m^{\times r})$ is non-special.
\end{proposition}

\begin{proof}
The proof is analogous. We begin with $m_0$-reduction
to obtain $(\overline{b},\{b\}^{\times \ell})$ for $b=d+1-m_0\geq a$
and some $\ell \in \N$. The last diagram can be reduced to some
$D \in \mathcal D$, or we end up with system without conditions.
\end{proof}

For a given $m$, $a$ and $k$ the set $\mathcal D$ defined
in Prop. \ref{downszer} or \ref{downup} can be very large.
On the other hand we do not need to consider diagrams, which
cannot be obtained as reductions of diagrams
of type $(\overline{a},\{a\}^{\times \ell})$.

\begin{proposition}
\label{throwout}
Let $D=(\dots,a,b,\dots)$ be a diagram obtained by
a sequence of $m$-reductions from a diagram $D=(\dots,a',b',\dots)$
($a'$ (resp. $b'$) stands on the same positions as $a$ (resp. $b$)).
If $b>0$ then
\begin{equation}
\tag{$\ast$}
\label{inthr}
a+(a-b+b'-a')m \geq a'. 
\end{equation}
\end{proposition}

\begin{proof}
Assume the contrary. Each reduction working on layer ``$a$''
works also on layer ``$b$'', moreover, the layer ``$b$''
is reduced stronger. Therefore at most $a-b+b'-a'$ such $m$-reductions
are possible, each one lowers the layer ``$a$'' by at most
$m$, which gives the size of this layer (at the beginning) at most
$a+(a-b+b'-a')m$, a contradiction.
\end{proof}

\begin{example}
Let
$\mathcal D = \{ (\overline{16},a_1,\dots,a_5) : 
16 \geq a_1 -1 \geq \dots \geq a_{5} -4 \}$.
We have $\# \mathcal D = 27896$, but if we throw out
all diagrams not satisfying inequality \eqref{inthr} we
will have $12799$ diagrams.
\end{example}

Observe that the dimension of $V(D;m_1,\dots,m_r)$ can be computed
by solving some (large) system of linear equations.
Usually this involves computation of the rank
of $\#D \times \sum_{j=1}^r \binom{m_j+1}{2}$ matrix (see e.g. \cite{redmd}).

We will use the following algorithm.

\begin{algorithm}{{\sc InitialCases}}

\noindent 
\begin{tabular}{rl}
{\bf Input:} & $m,a,k \in \N$ \\
{\bf Output:} & {\sc ok} or {\sc not ok}.
\end{tabular}
\\

\noindent
{\bf if} $a < m$ {\bf or} ($k=0$ {\bf and} $a<2m$) {\bf then} {\bf return} {\sc not ok};\\
{\bf if} $k=0$ {\bf then} prepare $\mathcal D = \{ (\overline{a},a_1,\dots,a_{m-1}) : 
a \geq a_1 -1 \geq a_2 -2 \geq \dots \geq a_{m-1} -(m-1) \}$;\\
{\bf if} $k>0$ {\bf then} prepare $\mathcal D = \{ (\overline{a},\{a\}^{\times k},a_1,\dots,a_{m-1}) : 
a \geq a_1 \geq a_2 \geq \dots \geq a_{m-1} \}$;\\
{\bf for each} $D \in \mathcal D$ satisfying inequality \eqref{inthr} {\bf do}\\
\hspace*{0.5cm} compute $p = p(D)$;\\
\hspace*{0.5cm} check non-specialty of $V(D;m^{\times p})$ and
$V(D;m^{(\times p+1)})$ by direct computation;\\
\hspace*{0.5cm} {\bf if} one of these systems is special {\bf then} {\bf return} {\sc not ok};\\
{\bf end for each}\\
{\bf return} {\sc ok};\\
\end{algorithm}

This algorithm has been implemented by the author (the work is done in
Free Pascal; the source code can be downloaded from \cite{MYWWW}).
In the table (see Tab. \ref{tab1}) we present the results of {\sc InitialCases} for $m=7$, $8$, $9$, $10$
and various $a$ and $k$. During implementation the following trick has been added.
Let $\mathcal D$ be a set of diagrams to decide about non-specialty, choose $s \geq 1$.
First, we compute
$$\mathcal D_{\red}^{s} = \{ \red_m^{s}(D) : D \in \mathcal D\},$$
where 
$$\red_m^s(D) \rdf \underbrace{\red_m(\red_m(\dots \red_m(D) \dots ))}_{s \text{ times}}.$$
Then, by matrix computation, create
$$\mathcal D_{\red, \text{ok}}^{s} = \{ D \in \mathcal D_{\red}^{s} : V(D;m^{\times p(D)}) \text{ and }
V(D;m^{\times (p(D)+1)}) \text{ are non-special}\}.$$
Let
$$\mathcal D_{\text{done}}^{s} = \{ D \in \mathcal D : \red_m^s(D) \in \mathcal D_{\red, \text{ok}}^{s} \}.$$
By reduction algorithm, $V(D;m^{\times p(D)})$ and $V(D;m^{\times (p(D)+1)})$
are non-special for all $D \in \mathcal D_{\text{done}}^{s}$.
We must check the remaining cases belonging to $\mathcal D \setminus \mathcal D_{\text{done}}^{s}$.
To do this proceed with new $\mathcal D \rdf \mathcal D \setminus \mathcal D_{\text{done}}^{s}$ and
new $s \rdf s-1$. For $s=0$ (final step) no reduction is performed, we only check non-specialty.

\begin{table}[ht!]
$$
\begin{array}{|c|c|c|c||c|c|c|c|c|c|}
\hline
m & a & k & \max_{D \in \mathcal D} p(D)+1 & m & a & k & \max_{D \in \mathcal D} p(D)+1 \\
\hline
7 & 18 & 0  &  11    &  9 & 20 & 4  &  11  \\
7 & 17 & 1  &  10    &  9 & 19 & 5  &  10  \\
7 & 16 & 2  &  10    &  9 & 18 & 6  &  10  \\
7 & 15 & 3  &  10    &  9 & 17 & 7  &  10  \\
7 & 14 & 4  &  9     &  9 & 16 & 8  &  9   \\
7 & 13 & 5  &  9     &  9 & 15 & 14 &  11  \\
7 & 12 & 11 &  11    &  9 & 14 & 17 &  11  \\
7 & 11 & 13 &  10    &  9 & 13 & 29 &  13  \\
7 & 10 & 24 &  13    &  9 & 12 & 62 &  21  \\
\hline
8 & 21 & 0  &  12    & 10 & 26 & 0  &  12  \\
8 & 20 & 1  &  11    & 10 & 25 & 1  &  11  \\
8 & 19 & 2  &  11    & 10 & 24 & 2  &  11  \\
8 & 18 & 3  &  10    & 10 & 23 & 3  &  11  \\
8 & 17 & 5  &  10    & 10 & 22 & 4  &  10  \\ 
8 & 16 & 5  &  10    & 10 & 21 & 5  &  10  \\
8 & 15 & 6  &  9     & 10 & 20 & 6  &  10  \\
8 & 14 & 7  &  9     & 10 & 19 & 7  &  9   \\
8 & 13 & 13 &  10    & 10 & 18 & 13 &  11  \\
8 & 12 & 19 &  11    & 10 & 17 & 15 &  11  \\
8 & 11 & 41 &  17    & 10 & 16 & 17 &  11  \\
\cline{1-4}
9 & 24 & 0  &  12    & 10 & 15 & 26 &  12  \\
9 & 23 & 1  &  11    & 10 & 14 & 41 &  15  \\
9 & 22 & 2  &  11    & 10 & 13 & 79 &  23  \\
\cline{5-8}
9 & 21 & 3  &  11    &    &    &    &      \\
\hline
\end{array}
$$
\caption{Results of {\sc InitialCases}. For all above the result was {\sc ok}}\label{tab1}
\end{table}

\section{Remaining cases}

According to Prop. \ref{downszer}, \ref{downup} and Tab. \ref{tab1}
we are left with the cases presented in Tab. \ref{tab2}.

\begin{table}[ht!]
$$
\begin{array}{|l|l||l|l|}
\hline
\text{system} & & \text{system} & \\
\hline
\hline
\sys(m+k;m,7^{\times r}) & \begin{array}{l} k \geq 17, \\ r \in \{9,10,11\} \end{array} & \sys(m+k;7^{\times r}) & \begin{array}{l} k \in \{16,15,14\}, \\ r \in \{9,10\} \end{array} \\
\hline
\sys(m+k;m,7^{\times 9}) & \begin{array}{l} k \in \{13,12\} \end{array} \begin{array}{l} \\ \\ \end{array} & \sys(m+11;m,7^{\times r}) & \begin{array}{l} r \in \{9,10,11\} \end{array} \\
\hline
\sys(m+10;m,7^{\times r}) & \begin{array}{l} r \in \{9,10\} \end{array} \begin{array}{l} \\ \\ \end{array} & \sys(m+9;m,7^{\times r}) & \begin{array}{l} r \in \{9,10,11,12,13\} \end{array} \\
\hline
\sys(m+k;m,7^{\times r}) & \begin{array}{l} k \in \{0,\dots,8\}, \\ r \geq 9 \end{array} & & \\
\hline
\hline
\sys(m+k;m,8^{\times r}) & \begin{array}{l} k \geq 20, \\ r \in \{9,\dots,12\} \end{array} & \sys(m+k;m,8^{\times r}) & \begin{array}{l} k \in \{17,16,15,12\}, \\ r \in \{9,10\} \end{array} \\
\hline
\sys(m+k;m,8^{\times r}) & \begin{array}{l} k \in \{19,18,11\}, \\ r \in \{9,10,11\} \end{array} & \sys(m+k;m,8^{\times 9}) & \begin{array}{l} k \in \{14,13\} \end{array} \\
\hline
\sys(m+10;m,8^{\times r}) & \begin{array}{l} r \in \{9,\dots,17\} \end{array} & \sys(m+k;m,8^{\times r}) & \begin{array}{l} k \in \{0,\dots,9\}, \\ r \geq 9 \end{array} \\
\hline
\hline
\sys(m+k;m,9^{\times r}) & \begin{array}{l} k \geq 23, \\ r \in \{9,\dots,12\} \end{array} & \sys(m+k;m,9^{\times r}) & \begin{array}{l} k \in \{22,\dots,19,14,13\}, \\ r \in \{9,10,11\} \end{array} \\
\hline
\sys(m+k;m,9^{\times r}) & \begin{array}{l} k \in \{18,17,16\}, \\ r \in \{9,10\} \end{array} & \sys(m+15;m,9^{\times 9}) & \\
\hline
\sys(m+12;m,9^{\times r}) & \begin{array}{l} r \in \{9,\dots,13\} \end{array} \begin{array}{l} \\ \\ \end{array} & \sys(m+11;m,9^{\times r}) & \begin{array}{l} r \in \{9,\dots,21\} \end{array} \\
\hline
\sys(m+k;m,9^{\times r}) & \begin{array}{l} k \in \{0,\dots,10\}, \\ r \geq 9 \end{array} & & \\
\hline
\hline
\sys(m+k;m,10^{\times r}) & \begin{array}{l} k \geq 25, \\ r \in \{9,\dots,12\} \end{array} & \sys(m+k;m,10^{\times r}) & \begin{array}{l} k \in \{24,23,22,17,16,15\}, \\ r \in \{9,10,11\} \end{array} \\
\hline
\sys(m+k;m,10^{\times r}) & \begin{array}{l} k \in \{21,20,19\}, \\ r \in \{9,10\} \end{array} & \sys(m+18;m,10^{\times 9}) & \\
\hline
\sys(m+14;m,10^{\times r}) & \begin{array}{l} r \in \{9,\dots,12\} \end{array} \begin{array}{l} \\ \\ \end{array} & \sys(m+13;m,10^{\times r}) & \begin{array}{l} r \in \{9,\dots,15\} \end{array} \\
\hline
\sys(m+12;m,10^{\times r}) & \begin{array}{l} r \in \{9,\dots,23\} \end{array} \begin{array}{l} \\ \\ \end{array} & \sys(m+k;m,10^{\times r}) & \begin{array}{l} k \in \{0,\dots,11\}, \\ r \geq 9 \end{array} \\
\hline
\end{array}
$$
\caption{Cases to be considered separately}\label{tab2}
\end{table}

In what follows we will solve all these cases. This may be boring; for every system
we must show that it is either non-special or $-1$-special.
We will use Cremona transformations, glueing (Thm. \ref{glue};
in most cases we glue four points) and
known facts about Harbourne-Hirschowitz conjecture
(such as it holds for multiplicities bounded by $11$). Observe that if a system
with non-negative multiplicities is in standard form 
and is based on at most $9$ points then it is non-special.

In what follows we write (for simplicity)
$\sys(d;\dots,m^{\times a,b,c,\dots})$ for a family of systems
$\{\sys(d;\dots,m^{\times \ell}) : \ell = a,b,c,\dots \}$.

The remaining cases can be divided with respect to methods of showing
non-specialty or $-1$-specialty. Therefore we present all used methods
(and an example for each of them), then, for each method, we give a list of cases
that can be done with this method.
For all considered systems we assume that $m \geq 12$.

\subsection{Glueing}
Glue four points $\sys(m+k;m,m_0^{\times r}) \to \sys(m+k;m,2m_0+1,m_0^{\times (r-4)})$.
The resulting system should be in standard form, based on at most $9$ points,
hence non-special, and with non-negative dimension.
As an example take $\sys(m+k;m,7^{\times 9,10,11})$,
$k \geq 22$, $m \geq 12$. After glueing we have
$\sys(m+k;m,15,7^{\times 5,6,7})$. This system is in standard form
since $m+k-m-15-7=m-22\geq 0$ and $m+k-15-7-7\geq 5$.
We also have $\vdim \sys(m+k;m,15,7^{\times 7}) = (k^2+2mk+3k+2m-630)/2-1 \geq 235$.
This method can be applied to the following systems:
$$
\begin{array}{ll@{\qquad\qquad}ll}
\sys(m+k;m,7^{\times 9,10,11}), & k \geq 22, &
\sys(m+k;m,8^{\times 9,10,11}), & k \geq 25, \\
\sys(m+k;m,9^{\times 9,10,11}), & k \geq 28, &
\sys(m+k;m,10^{\times 9,10,11}), & k \geq 31.
\end{array}
$$

\subsection{Double glueing}
As before, but we must glue twice (i.e. $8$ points
of multiplicity $m_0$ to $2$ points of multiplicity $2m_0+1$).
As an example consider
$\sys(m+k;m,8^{\times 12})$ for $k \geq 34$, $m \geq 12$.
Glue twice to obtain $\sys(m+k;m,17^{\times 2},8^{\times 4})$ in standard form.
This method can be applied to the following systems:
$$
\begin{array}{ll@{\qquad\qquad}ll}
\sys(m+k;m,8^{\times 12}), & k \geq 34, &
\sys(m+k;m,9^{\times 12}), & k \geq 38, \\
\sys(m+k;m,10^{\times 12}), & k \geq 42.
\end{array}
$$

\subsection{Glueing and Cremona}
Glue four points $\sys(m+k;m,m_0^{\times r}) \to \sys(m+k;m,2m_0+1,m_0^{\times (r-4)})$.
Then use Cremona transformation based on points
with multiplicities $m$, $2m_0+1$ and $m_0$ obtaining
$\sys(m+2k-3m_0-1;m+k-3m_0-1,k-m_0,k-2m_0-1,m_0^{\times (r-5)})$.
The last system should be non-special in standard form.
As an example take $\sys(m+k;m,7^{\times 9,10,11})$,
$k \in \{17,\dots,21\}$, $m \geq 12$. After glueing we have
$\Cr(\sys(m+k;m,15,7^{\times 5,6,7}))=\sys(m+2k-22;m+k-22,k-7,k-15,7^{\times 4,5,6})$.
This system is in standard form since
$(m+2k-22)-(m+k-22)-(k-7)-(k-15)=-k+22 \geq 0$,
$(m+2k-22)-(m+k-22)-(k-7)-7=0$ and $(m+2k-22)-(k-7)-7-7=m+k-29\geq 0$.
The computation of virtual dimension is straightforward and gives
$\vdim = (k^2+2mk+3k+2m-630)/2-1 \geq 70$.
This method can be applied to the following systems:
$$
\begin{array}{ll@{\qquad\qquad}ll}
\sys(m+k;m,7^{\times 9,10,11}), & k \in \{17,\dots,21\}, &
\sys(m+k;m,8^{\times 9,10,11}), & k \in \{20,\dots,24\}, \\
\sys(m+k;m,9^{\times 9,10,11}), & k \in \{25,26,27\}, & 
\sys(m+k;m,10^{\times 9,10,11}), & k \in \{28,29,30\}.
\end{array}
$$

\subsection{Double glueing and Cremona}
As before, but we must glue twice.
As an example consider
$\sys(m+k;m,8^{\times 12})$ for $k \in \{23,\dots,33\}$, $m \geq 12$.
Glue twice to obtain 
$\Cr(\sys(m+k;m,17^{\times 2},8^{\times 4}))=\sys(m+2k-34;m+k-34,(k-17)^{\times 2},8^{\times 4})$
in standard form.
This method can be applied to the following systems:
$$
\begin{array}{ll@{\qquad\qquad}ll}
\sys(m+k;m,8^{\times 12}), & k \in \{23,\dots,33\}, &
\sys(m+k;m,9^{\times 12}), & k \in \{25,\dots,37\}, \\
\sys(m+k;m,10^{\times 12}), & k \in \{30,\dots,41\}.
\end{array}
$$

\subsection{Glueing and Cremona(s)}
Glue four points $\sys(m+k;m,m_0^{\times r}) \to \sys(m+k;m,2m_0+1,m_0^{\times (r-4)})$.
Then use Cremona transformation based on points
with multiplicities $m$, $2m_0+1$ and $m_0$ obtaining
$\sys(m+2k-3m_0-1;m+k-3m_0-1,k-m_0,k-2m_0-1,m_0^{\times (r-5)})$.
The last system should be non-special in standard form except for a finite
number of cases for low values of $m$. For each of these cases we must use
an additional sequence of Cremona transformations to end up with a system
in standard form.
As an example take $\sys(m+16;m,7^{\times 9,10})$
for $m \geq 12$. After glueing we have
$\Cr(\sys(m+16;m,15,7^{\times 5,6}))=\sys(m+10;m-6,9,7^{\times 4,5},1)$.
Since $m+10-9-7-7=m-13$ the last system is in standard form for $m\geq 13$.
The remaining case $m=12$ can be done as follows:
$\CrSt(\sys(22;6,9,7^{\times 4,5},1))=\sys(20;7^{\times 1,2},6^{\times 5},1)$.
This method can be applied to the following systems:
$$
\begin{array}{ll@{\qquad\qquad}ll}
\sys(m+k;m,7^{\times 9,10}), & k \in \{16,15,14\}, &
\sys(m+k;m,8^{\times 9,10,11}), & k \in \{19,18\}, \\
\sys(m+17;m,8^{\times 9,10}), & &
\sys(m+16;m,8^{\times 9,10}), & m \geq 13, \\
\sys(m+k;m,9^{\times 9,10,11}), & k \in \{24,23,22,21\}, &
\sys(m+20;m,9^{\times 9,10}), &  \\
\sys(m+20;m,9^{\times 11}), & m \geq 14, &
\sys(m+19;m,9^{\times 9,10}), & m \geq 13, \\
\sys(m+19;m,9^{\times 11}), & m \geq 15, &
\sys(m+18;m,9^{\times 9,10}), & m \geq 14, \\
\sys(m+k;m,10^{\times 9,10,11}), & k \in \{25,26,27\}, &
\sys(m+24;m,10^{\times 9,10,11}), & \\
\sys(m+23;m,10^{\times 9,10}), & &
\sys(m+23;m,10^{\times 11}), & m \geq 14, \\
\sys(m+22;m,10^{\times 9,10}), & m \geq 13, &
\sys(m+22;m,10^{\times 11}), & m \geq 15, \\
\sys(m+21;m,10^{\times 9,10}), & m \geq 14, &
\sys(m+20;m,10^{\times 9,10}), & m \geq 16.
\end{array}
$$

\subsection{Double glueing and Cremona(s)}
As before, but we must glue twice.
As an example consider
$\sys(m+22;m,8^{\times 12})$ for $m \geq 12$.
Glue twice to obtain 
$\Cr(\sys(m+22;m,17^{\times 2},8^{\times 4}))=\sys(m+10;m-12,8^{\times 4},5^{\times 2})$.
For $m \geq 14$ the last system is in standard form, the remaining cases are
$$
\begin{array}{ll}
\CrSt(\sys(23;8^{\times 4},5^{\times 2},1))=\sys(22;8,7^{\times 3},5^{\times 2},1),&
\CrSt(\sys(22;8^{\times 4},5^{\times 2}))=\sys(20;8,6^{\times 3},5^{\times 2}).
\end{array}
$$
This method can be applied to the following systems:
$$
\begin{array}{ll@{\qquad\qquad}ll}
\sys(m+22;m,8^{\times 12}), &  &
\sys(m+21;m,8^{\times 12}), &  \\
\sys(m+20;m,8^{\times 12}), & m \geq 13, &
\sys(m+24;m,9^{\times 12}), & m \geq 13, \\
\sys(m+23;m,9^{\times 12}), & m \geq 14, &
\sys(m+29;m,10^{\times 12}), & \\
\sys(m+28;m,10^{\times 12}), & m \geq 13, &
\sys(m+27;m,10^{\times 12}), & m \geq 14, \\
\sys(m+26;m,10^{\times 12}), & m \geq 15, &
\sys(m+25;m,10^{\times 12}), & m \geq 16.
\end{array}
$$

\subsection{Glue, Cremona, glue, Cremona}
Glue four points of equal multiplicity, then perform
Cremona transformation several times to obtain system with lower
multiplicities. Then glue four points (but now the multiplicities
are lower) and use Cremona transformation(s) to obtain a non-special
system in standard form.
As an example consider
$\sys(32;12,8^{\times 12})$.
Glue to obtain
$\CrSt(\sys(32;17,12,8^{\times 8}))=\sys(24;9,8,7^{\times 7},3)$.
Glue again to consider
$\CrSt(\sys(24;15,9,8,7^{\times 3},3))=\sys(10;6,2^{\times 3},1^{\times 2})$.
This method can be applied to the following systems:
$$
\begin{array}{ll@{\qquad\qquad}ll}
\sys(32;12,8^{\times 12}), & &
\sys(36;12,9^{\times 12}), & \\
\sys(36;13,9^{\times 12}), & &
\sys(35;12,9^{\times 12}), & \\
\sys(40;12,10^{\times 12}), & &
\sys(40;13,10^{\times 12}), & \\
\sys(39;12,10^{\times 12}), & &
\sys(40;14,10^{\times 12}), & \\
\sys(39;13,10^{\times 12}), & &
\sys(38;12,10^{\times 12}), & \\
\sys(40;15,10^{\times 12}), & &
\sys(39;14,10^{\times 12}), & \\
\sys(38;13,10^{\times 12}). & &
\end{array}
$$

\subsection{Cremona (even) and glueing}
Let us consider $\sys(m+k;m,m_0^{\times 2r})$ such that $k-2m_0<0$.
Perform Cremona transformations based on the first point and two points
with multiplicity $m_0$. Each time the degree and the first 
multiplicity is changed by $k-2m_0$. We end up with
$\sys(m+k+r(k-2m_0);m+r(k-2m_0),(k-m_0)^{\times 2r})$.
For $m$ such that $m+r(k-2m_0)\leq 11$ the situation is known
(observe that this multiplicity can be negative).
Otherwise glue four points of multiplicity $k-m_0$ and end up with
system in standard form.
As an example consider $\sys(m+9;m,7^{\times 10})$ for $m \geq 12$.
Use Cremona transformations to obtain $\sys(m-16;m-25,2^{\times 10})$.
For $m \geq 37$ glue points to obtain
$\sys(m-16;m-25,5,2^{\times 6})$ in standard form.
This method can be applied to the following systems:
$$
\begin{array}{ll@{\qquad\qquad}ll}
\sys(m+10;m,7^{\times 10}), & &
\sys(m+9;m,7^{\times 10}), & \\
\sys(m+11;m,8^{\times 10}), & &
\sys(m+10;m,8^{\times 10}), & \\
\sys(m+13;m,9^{\times 10}), & &
\sys(m+12;m,9^{\times 10}), & \\
\sys(m+14;m,10^{\times 10}), & &
\sys(m+13;m,10^{\times 10}), & \\
\sys(m+12;m,10^{\times 10}). & &
\end{array}
$$

\subsection{Cremona (even) and multiple glueing}
As before, but we must glue several times to produce system
based on at most $9$ points.
As an example consider
$\sys(m+9;m,7^{\times 12})$ for $m \geq 12$.
Use Cremona transformations
to obtain $\sys(m-21;m-30,2^{\times 12})$. For $m \geq 42$ glue twice and
finish with $\sys(m-21;m-30,5^{\times 2},2^{\times 4})$ in standard form.
This method can be applied to the following systems (in square
brackets we indicate how many times we glue):
$$
\begin{array}{ll@{\qquad\qquad}ll}
\sys(m+9;m,7^{\times 12}), & [2], &
\sys(m+10;m,8^{\times 12,14}), & [2], \\
\sys(m+11;m,9^{\times 10,12,14}), & [2], &
\sys(m+11;m,9^{\times 16,18,20}), & [4], \\
\sys(m+12;m,10^{\times 12,14}), & [2], &
\sys(m+12;m,10^{\times 16,18}), & [4], \\
\sys(m+12;m,10^{\times 20,22}), & [5]. &
\end{array}
$$

\subsection{Cremona (even), glueing and Cremona(s)}
As before, consider 
$\sys(m+k;m,m_0^{\times 2r})$ such that $k-2m_0<0$,
but now the system after glueing
$\sys(m+k+r(k-2m_0);m+r(k-2m_0),2k-2m_0+1,(k-m_0)^{\times (2r-4)})$ is not
in standard form, we must use another Cremona transformation
based on points with multiplicities
$m+r(k-2m_0)$, $2k-2m_0+1$, $k-m_0$.
Now the system is in standard form for $m$ big enough, for a finite
set of values of $m$ (this set may be empty) we must use
additional sequence of Cremona transformation(s).
As an example consider
$\sys(m+11;m,7^{\times 10})$ for $m \geq 12$.
Transform this system
into $\sys(m-4;m-15,4^{\times 10})$. For $m \geq 27$
glue points to $\sys(m-4;m-15,9,4^{\times 6})$. The standard form
is $\sys(m-6;m-17,7,4^{\times 5},2)$.
Another example is $\sys(m+15;m,8^{\times 10})$.
This system can be transformed
into $\sys(m+10;m-5,7^{\times 10})$. For $m \geq 17$ use glueing to
obtain $\sys(m+10;m-5,15,7^{\times 6})$. For $m \geq 19$
the standard form is $\sys(m+3;m-12,8,7^{\times 5})$, the remaining cases are
$$
\begin{array}{ll}
\CrSt(\sys(20;8,7^{\times 5},5))=\sys(14;5^{\times 2},4^{\times 5}), &
\CrSt(\sys(21;8,7^{\times 5},6))=\sys(19;7,6^{\times 6}).
\end{array}
$$
This method can be applied to the following systems:
$$
\begin{array}{ll@{\qquad\qquad}ll}
\sys(m+11;m,7^{\times 10}), & &
\sys(m+15;m,8^{\times 10}), & \\
\sys(m+12;m,8^{\times 10}), & &
\sys(m+17;m,9^{\times 10}), & \\
\sys(m+16;m,9^{\times 10}), & &
\sys(m+14;m,9^{\times 10}), & \\
\sys(m+19;m,10^{\times 10}), & m \neq 17, &
\sys(m+17;m,10^{\times 10}), & \\
\sys(m+16;m,10^{\times 10}), & &
\sys(m+15;m,10^{\times 10}). & \\
\end{array}
$$

\subsection{Cremona (even), multiple glueing and Cremona(s)}
As before, but we must glue several times.
As an example consider
$\sys(m+10;m,8^{\times 16})$ for $m \geq 12$.
This system can be transformed
into $\sys(m-38;m-48,2^{\times 16})$. For $m \geq 60$ glue three times
to obtain $\sys(m-38;m-48,5^{\times 3},2^{\times 4})$, which can be transformed into
the standard form $\sys(m-43;m-53,2^{\times 4})$.
This method can be applied to the following systems (in square
brackets we indicate how many times we glue):
$$
\begin{array}{ll@{\qquad\qquad}ll}
\sys(m+10;m,8^{\times 16}), & [3], &
\sys(m+12;m,9^{\times 12}), & [2], \\
\sys(m+14;m,10^{\times 12}), & [2], &
\sys(m+13;m,10^{\times 12,14}), & [2]. \\
\end{array}
$$

\subsection{Cremona (odd), glueing and Cremona(s)}
Consider 
$\sys(m+k;m,m_0^{\times 2r+1})$ such that $k-2m_0<0$. This system
can be transformed into 
$\sys(m+k+r(k-2m_0);m+r(k-2m_0),m_0,(k-m_0)^{\times 2r})$.
For $m$ such that $m+r(k-2m_0)\leq 11$ the situation is known.
Otherwise glue four points of multiplicity $k-m_0$ to obtain
$\sys(m+k+r(k-2m_0);m+r(k-2m_0),2k-2m_0+1,m_0,(k-m_0)^{\times (2r-4)})$.
Use another Cremona transformation
based on points with multiplicities
$m+r(k-2m_0)$, $2k-2m_0+1$, $m_0$ to obtain the system
$\sys(m+r(k-2m_0)+m_0-1;m+m_0-k+r(k-2m_0)-1,k-m_0,2m_0-k-1,(k-m_0)^{\times (2r-4)})$
in standard form for $m$ big enough.
For remaining values of $m$ we must use additional Cremona
transformation(s) to end up in standard form.
As an example consider 
$\sys(m+13;m,7^{\times 9})$ for $m \geq 12$.
This system can be transformed
into $\sys(m+9;m-4,7,6^{\times 8})$. For $m < 16$ the situation is known.
For $m \geq 16$ let us glue points to
$\sys(m+9;m-4,13,7,6^{\times 4})$. The standard form of the last system
is $\sys(m+2;m-11,6^{\times 5})$.
This method can be applied to the following systems:
$$
\begin{array}{ll@{\qquad\qquad}ll}
\sys(m+k;m,7^{\times 9}), & k \in \{13,12,10\}, &
\sys(m+k;m,7^{\times 9,11}), & k \in \{11,9\}, \\
\sys(m+k;m,8^{\times 9}), & k \in \{12,\dots,15\}, &
\sys(m+k;m,8^{\times 9,11}), & k \in \{11,10\}, \\
\sys(m+k;m,9^{\times 9}), & k \in \{17,16,15\}, &
\sys(m+k;m,9^{\times 9,11}), & k \in \{14,13,12\}, \\
\sys(m+19;m,10^{\times 9}), & m \neq 16, &
\sys(m+18;m,10^{\times 9}), & \\
\sys(m+k;m,10^{\times 9,11}), & k \in \{12,\dots,17\}. &
\end{array}
$$

\subsection{Cremona (odd), multiple glueing and Cremona(s)}
As before, but we must glue several times to obtain the system with
at most $9$ multiplicities.
As an example consider
$\sys(m+9;m,7^{\times 13})$ for $m \geq 12$.
Transform our system into $\sys(m-21;m-30,7,2^{\times 12})$. For $m \geq 42$
glue twice to obtain
$\Cr(\sys(m-21;m-30,7,5,5,2^{\times 4}))=\sys(m-24;m-33,5,4,2^{\times 5})$
in standard form.
This method can be applied to the following systems (in square
brackets we indicate how many times we glue):
$$
\begin{array}{ll@{\qquad\qquad}ll}
\sys(m+9;m,7^{\times 13}), & [2], &
\sys(m+10;m,8^{\times 13}), & [2], \\
\sys(m+10;m,8^{\times 15,17}), & [3], &
\sys(m+12;m,9^{\times 13}), & [2], \\
\sys(m+11;m,9^{\times 9,11,13}), & [2], &
\sys(m+11;m,9^{\times 15}), & [3], \\
\sys(m+11;m,9^{\times 17,19}), & [4], &
\sys(m+11;m,9^{\times 21}), & [5], \\
\sys(m+13;m,10^{\times 13}), & [2], &
\sys(m+13;m,10^{\times 15}), & [3], \\
\sys(m+12;m,10^{\times 13}), & [2], &
\sys(m+12;m,10^{\times 15,17}), & [3], \\
\sys(m+12;m,10^{\times 19}), & [4], &
\sys(m+12;m,10^{\times 21,23}), & [5]. \\
\end{array}
$$

\subsection{Negative glueing and Cremona}
We use glueing to show that the system
$\sys(m+k;m,m_0^{\times r})$ with negative virtual dimension is empty.
Therefore we glue four points of multiplicity $m_0$ to
one point of multiplicity $2m_0$, then we use Cremona
transformation(s) to show that the resulting system is empty.
As an example consider
$\sys(32;13,9^{\times 11})$.
Use glueing to consider $\sys(32;18,13,9^{\times 7})$ which can be tranformed
into $\sys(0;2^{\times 3},1,(-1)^{\times 5},-2,-4)$. The last system is empty.
This method can be applied to the following systems:
$$
\begin{array}{ll@{\qquad\qquad}ll}
\sys(32;12,9^{\times 11}), & &
\sys(32;13,9^{\times 11}), & \\
\sys(31;12,9^{\times 11}), & &
\sys(31;13,9^{\times 10}), & \\
\sys(30;12,9^{\times 10}), & &
\sys(35;12,10^{\times 11}), & \\
\sys(35;13,10^{\times 11}), & &
\sys(34;12,10^{\times 11}), & \\
\sys(34;13,10^{\times 10}), & &
\sys(33;12,10^{\times 10}), & \\
\sys(34;14,10^{\times 10}), & &
\sys(33;13,10^{\times 10}), & \\
\sys(32;12,10^{\times 10}). & &
\end{array}
$$

\subsection{Low multiplicities}
Consider $\sys(m+k;m,m_0^{\times r})$ for $k-m_0 \leq 1$.
As before, this system can be transformed (by a sequence of Cremona transformations)
to a system in standard form with at most two arbitrary ``big'' multiplicities, the other being
strictly less that $2$. Let $L=\sys(d;m_1,m_2,m_0^{\times s})$ be such a system,
$m_0 \leq 1$. If $m_0 \leq -2$ then $L$ is $-1$-special if and only
if $L$ is non-empty, which is equivalent to non-emptyness of $\sys(d;m_1,m_2)$.
For $m_0=-1$ or $m_0=0$ it is enough to
consider $\sys(d;m_1,m_2)$ based on at most two points. For $m_0=1$ we have
two cases. If $m_1 \geq -1$, $m_2 \geq -1$ then $L$ is non-special
since multiplicity $1$ always imposes an independent condition.
For the opposite case we must decide whether $L$ is non-empty. Dropping
negative multiplicities we end up with a system with at most one
multiplicity not equal to $1$.
As an example consider
$\sys(m+8;m,7^{\times 2r+1})$.
This system can be transformed to $\sys(m+k-6r;m-6r,7,1^{\times 2r})$
in standard form.
If $m-6r \geq -1$ then the system is non-special.
If $m-6r < -1$ then the system is $-1$ special if and only if
$\sys(m+8-6r;7,1^{\times 2r})$ is non-empty, which holds for
$\binom{m-6r+10}{2} \geq 28+2r$. In fact we have
$\binom{m-6r+10}{2} \leq \binom{8}{2} = 28$, so $r=0$ and our
system is non-special.
This method can be applied to the following systems:
$$
\begin{array}{ll@{\qquad\qquad}ll}
\sys(m+k;m,7^{\times r}), & k \leq 8, \, r \geq 9, & 
\sys(m+k;m,8^{\times r}), & k \leq 9, \, r \geq 9, \\
\sys(m+k;m,9^{\times r}), & k \leq 10, \, r \geq 9, &
\sys(m+k;m,10^{\times r}), & k \leq 11, \, r \geq 9.
\end{array}
$$

\subsection{Additional methods}
We use non-standard glueing, reduction algorithm, etc.

\paragraph{$\sys(28;12,8^{\times 9})$}
Glue three points
(using non-special system $\sys(15;8^{\times 3})$) to obtain
$\sys(28;16,12,8^{\times 6})$. The standard form of the last
system is $\sys(4;4)$.
\epar

\paragraph{$\sys(31;12,9^{\times 9})$}
This system is non-special due to reduction algorithm. We begin with
the diagram $(\overline{32})$, use $12$-reduction followed by
nine $9$-reductions. The last diagram is equal to $(\overline{6},6,6,5,5,2)$.
\epar

\paragraph{$\sys(31;13,9^{\times 9})$}
It is enough to show that $\sys(30;13,9^{\times 9})$ is non-special
(observe that the last system has virtual dimension equal to $-1$).
We have $\CrSt(\sys(30;13,9^{\times 9}))=\sys(26;9^{\times 2},8^{\times 8})$.
Since all the multiplicities are bounded by $11$ we can use Thm. \ref{up11}.
\epar

\paragraph{$\sys(34;12,10^{\times 9})$}
This system has positive virtual dimension. Since
$\sys(19;10^{\times 3})$ is non-empty and non-special we use glueing
to consider $\CrSt(\sys(34;20,12,10^{\times 6}))=\sys(4;2)$.
\epar

\paragraph{$\sys(35;15,10^{\times 9})$}
This system has positive virtual dimension. Since
$\sys(19;10^{\times 3})$ is non-empty and non-special we use glueing
to consider $\CrSt(\sys(35;20,15,10^{\times 6}))=\sys(5;5)$.
\epar

\paragraph{$\sys(32;12,10^{\times 9})$}
This system is empty due to reduction algorithm. We begin with
the diagram $(\overline{33})$, use $12$-reduction followed by
eight $10$-reductions. The last diagram is equal to $(\overline{6},6,6,5,5)$,
which can be enlarged to $(\overline{10})$ and reduced to an empty diagram.
\epar

\paragraph{$\sys(35;16,10^{\times 9})$}
This system can be transformed into $\sys(31;12,10,9^{\times 8})$. It is
enough to show that the system $\sys(30;12,10,9^{\times 8})$ is non-empty and non-special.
The last system can be transformed into $\sys(29;11,9^{\times 8},8)$ which
is non-special due to Thm. \ref{up11}.
\epar

\subsection{Direct computations}
Sometimes we are forced to compute the rank of the matrix associated
to a system. To make this task possible, we specialize to
random points and compute over $\mathbb F_p$ for some prime $p$.
If the rank is maximal for specialized points over $\mathbb F_p$ then it is maximal
over $\mathbb Q$ (and hence over any field of characteristic zero)
and for points in general position. Alternatively, we may use
diagram cutting method presented in \cite{md} (the author checked that
in all cases it is possible).
This method must be applied to the following systems:
$$
\begin{array}{ll@{\qquad\qquad}ll}
\sys(28;12,8^{\times 10}), & &
\sys(33;13,9^{\times 11}), & \\
\sys(33;14,9^{\times 11}), & &
\sys(31;12,9^{\times 10}), & \\
\sys(32;14,9^{\times 10}), & &
\sys(30;12,9^{\times 9}), & \\
\sys(37;12,10^{\times 12}), & &
\sys(36,13,10^{\times 11}), & \\
\sys(36;14,10^{\times 11}), & &
\sys(34;12,10^{\times 10}), & \\
\sys(34;13,10^{\times 9}), & &
\sys(33;12,10^{\times 9}), & \\
\sys(35;15,10^{\times 10}), & &
\sys(34;14,10^{\times 9}), & \\
\sys(33;13,10^{\times 9}), & &
\sys(36;17,10^{\times 10}). & \\
\end{array}
$$

\section{A note on Seibert's proof for $m=4$}

In \cite{sei4} all special systems of the form $\sys(d;m,4^{\times r})$ has
been classified. For all non-special cases but one the proof involved
techniques avoiding computation of the rank of matrix.
For $\sys(13;5,4^{\times 9})$ the author of \cite{sei4} used Maple program to compute
the rank of $105 \times 105$ matrix. The rank appeared to be maximal,
so the system is non-special. Using diagram cutting method
(introduced in \cite{mdwj}) we propose much nicer proof of this fact,
which can be easily checked by hand. The cutting is presented on
Fig. \ref{figcut}. The order of cutting is indicated by numbers
on a diagram. Now, Seibert's proof do not rely on
computation that cannot be done by hand.

\begin{figure}[ht!]
$$
\begin{texdraw}
\drawdim pt
\textref h:C v:C
\move(0 0)
\lvec(0 140)
\move(10 0)
\lvec(10 140)
\move(20 0)
\lvec(20 130)
\move(30 0)
\lvec(30 120)
\move(40 0)
\lvec(40 110)
\move(50 0)
\lvec(50 100)
\move(60 0)
\lvec(60 90)
\move(70 0)
\lvec(70 80)
\move(80 0)
\lvec(80 70)
\move(90 0)
\lvec(90 60)
\move(100 0)
\lvec(100 50)
\move(110 0)
\lvec(110 40)
\move(120 0)
\lvec(120 30)
\move(130 0)
\lvec(130 20)
\move(140 0)
\lvec(140 10)

\move(0 140)
\lvec(10 140)
\move(0 130)
\lvec(20 130)
\move(0 120)
\lvec(30 120)
\move(0 110)
\lvec(40 110)
\move(0 100)
\lvec(50 100)
\move(0 90)
\lvec(60 90)
\move(0 80)
\lvec(70 80)
\move(0 70)
\lvec(80 70)
\move(0 60)
\lvec(90 60)
\move(0 50)
\lvec(100 50)
\move(0 40)
\lvec(110 40)
\move(0 30)
\lvec(120 30)
\move(0 20)
\lvec(130 20)
\move(0 10)
\lvec(140 10)
\move(0 0)
\lvec(140 0)

\htext(5 135){2}
\htext(5 125){2}
\htext(5 115){2}
\htext(5 105){2}
\htext(5 95){6}
\htext(5 85){6}
\htext(5 75){6}
\htext(5 65){6}
\htext(5 55){8}
\htext(5 45){1}
\htext(5 35){1}
\htext(5 25){1}
\htext(5 15){1}
\htext(5 5){1}
\htext(15 125){2}
\htext(15 115){2}
\htext(15 105){2}
\htext(15 95){6}
\htext(15 85){6}
\htext(15 75){6}
\htext(15 65){8}
\htext(15 55){8}
\htext(15 45){8}
\htext(15 35){1}
\htext(15 25){1}
\htext(15 15){1}
\htext(15 5){1}
\htext(25 115){2}
\htext(25 105){2}
\htext(25 95){6}
\htext(25 85){6}
\htext(25 75){8}
\htext(25 65){8}
\htext(25 55){8}
\htext(25 45){8}
\htext(25 35){4}
\htext(25 25){1}
\htext(25 15){1}
\htext(25 5){1}
\htext(35 105){2}
\htext(35 95){6}
\htext(35 65){9}
\htext(35 55){8}
\htext(35 45){8}
\htext(35 35){7}
\htext(35 25){4}
\htext(35 15){1}
\htext(35 5){1}
\htext(45 55){9}
\htext(45 45){9}
\htext(45 35){7}
\htext(45 25){4}
\htext(45 15){4}
\htext(45 5){1}
\htext(55 55){9}
\htext(55 45){9}
\htext(55 35){7}
\htext(55 25){7}
\htext(55 15){4}
\htext(55 5){4}
\htext(65 75){5}
\htext(65 55){9}
\htext(65 45){9}
\htext(65 35){9}
\htext(65 25){7}
\htext(65 15){7}
\htext(65 5){4}
\htext(75 65){5}
\htext(75 55){5}
\htext(75 45){9}
\htext(75 35){9}
\htext(75 25){7}
\htext(75 15){7}
\htext(75 5){4}
\htext(85 55){5}
\htext(85 45){5}
\htext(85 35){5}
\htext(85 25){7}
\htext(85 15){7}
\htext(85 5){4}
\htext(95 45){5}
\htext(95 35){5}
\htext(95 25){5}
\htext(95 15){5}
\htext(95 5){7}
\htext(105 35){3}
\htext(105 25){3}
\htext(105 15){3}
\htext(105 5){3}
\htext(115 25){3}
\htext(115 15){3}
\htext(115 5){3}
\htext(125 15){3}
\htext(125 5){3}
\htext(135 5){3}
\end{texdraw}
$$
\caption{Diagram cutting for $\sys(13;5,4^{\times 9})$}\label{figcut}
\end{figure}

\end{document}